\documentclass[a4paper,draft,12pt,reqno]{article}
\usepackage{amsmath,amssymb,amsthm,xspace,amscd,latexsym,amssymb}
\usepackage[all]{xy}
\usepackage{enumerate}

\theoremstyle{plain}
\newtheorem{theorem}{Theorem}
\newtheorem{lemma}[theorem]{Lemma}
\newtheorem{proposition}[theorem]{Proposition}
\newtheorem{corollary}[theorem]{Corollary}

\newcommand{\tto}{\twoheadrightarrow}

\font\cs=rsfs10 at 12pt
\newcommand{\mcat}[1]{\mbox{\cs #1}\hspace{1.5pt}}
\newcommand{\LL}{\mcat L}

\DeclareMathOperator{\Hom}{Hom}

\begin{document} 

\title{A twisted approach to Kostant's problem}
\author{Volodymyr Mazorchuk}
\date{}
\maketitle

\begin{abstract}
We use Arkhipov's twisting functors to show that the universal enveloping algebra of
a semi-simple complex finite-dimensional Lie algebra surjects onto the space of
ad-finite endomorphisms of the simple highest weight module $L(\lambda)$, whose
highest weight is associated (in the natural way) with a subset of simple roots
and a simple root in this subset. This is a new step towards a complete answer 
to a classical question of Kostant. We also show how one can use the twisting 
functors to reprove the classical results related to this question.
\end{abstract}

\section{Introduction and notation}\label{s1}

Let $\mathfrak{g}$ be a complex semi-simple finite-dimensional Lie algebra with a
fixed triangular decomposition, 
$\mathfrak{g}=\mathfrak{n}_-\oplus \mathfrak{h}\oplus\mathfrak{n}_+$, and 
$U(\mathfrak{g})$ be its universal enveloping algebra. Then for every
two $\mathfrak{g}$-modules $M$ and $N$ the space $\Hom_{\mathbb{C}}(M,N)$ can
be viewed as a $U(\mathfrak{g})$-bimodule in the natural way (with the right
action of $U(\mathfrak{g})$  defined via the Chevalley involution). This bimodule then 
also becomes a $\mathfrak{g}$-module under the adjoint action. The bimodule 
$\Hom_{\mathbb{C}}(M,N)$ has a sub-bimodule, usually denoted by $\LL(M,N)$
(see for example \cite[Kapitel~6]{Ja}), 
which consists of all elements, the adjoint action of $U(\mathfrak{g})$ on which
is locally finite. Since $U(\mathfrak{g})$ itself consists of locally finite
elements under the adjoint action, it naturally maps to $\LL(M,M)$
for every $\mathfrak{g}$-module $M$, and the kernel of this map is the annihilator
$\mathrm{Ann}(M)$ of $M$ in $U(\mathfrak{g})$. The classical problem of Kostant (see 
for example \cite{Jo2}) is formulated  in the following way: 
\vspace{4mm}

{\em For which simple $\mathfrak{g}$-modules $M$ the natural injection
\begin{displaymath}
U(\mathfrak{g})/\mathrm{Ann}(M)\hookrightarrow \LL(M,M)
\end{displaymath}
is surjective?}
\vspace{4mm}

The complete answer to this problem is not known even for simple highest weight
modules. However, it is known that there are simple highest weight modules for which 
the answer is negative (see for example \cite[9.5]{Jo2}). There is also a  
classical class of simple highest weight modules, for which the answer is positive. 
It consists of all simple highest weight modules, whose highest weights are obtained 
from the antidominant one, applying the longest element of some parabolic subgroup 
of the Weyl group, see  \cite{GJ,Jo2,Ja}.

In the present paper we propose an approach to this problem, which uses Arkhipov's
twisting functors, see \cite{Ar}, and is based on the properties of these functors
obtained in \cite{AS}. In \cite{KM} it was shown that Arkhipov's functors are 
adjoint to Joseph's completion functors, see \cite{Jo1}, which suggests a close
connection to Kostant's problem. We base our arguments mostly on the results of 
\cite{AS} and also use some results from \cite{Kh,KM,MS}. The main properties of the 
twisting functors which we use are: the combinatorics of their action on Verma 
modules and the fact that they define a self-equivalence of the bounded derived 
category $\mathcal{D}^b(\mathcal{O})$ of the BGG-category 
$\mathcal{O}$. All this can be found in \cite{AS}.

Let $R$ be the root system of $\mathfrak{g}$ with the basis $B$, which corresponds to
the triangular decomposition above. Let further $W$ denote the Weyl group of 
$\mathfrak{g}$ with the identity element $e$. 
Then $W$ acts on $\mathfrak{h}^*$ both in the natural way (i.e. $\lambda\mapsto w(\lambda)$ 
for $\lambda\in  \mathfrak{h}^*$ and $w\in W$) and via the {\em dot action} defined as 
follows: $w\cdot\lambda=w(\lambda+\rho)-\rho$, $\lambda\in  \mathfrak{h}^*$, $w\in W$, 
where $\rho$ is the half of the sum of all positive roots. For $\alpha\in R$
denote by $s_{\alpha}$ the corresponding simple reflection, and for a reflection,
$s\in W$, we let $\alpha_s\in R$ be such that $s=s_{\alpha_s}$.
Fix some Weyl-Chevalley basis in $\mathfrak{g}$, say 
\begin{displaymath}
\{X_{\alpha}\,:\,\alpha\in R\}\cup \{H_{\beta}\,:\,\beta\in B\},
\end{displaymath}
and define $H_{\alpha}$, $\alpha\in R$, in the usual way.

For $\lambda\in \mathfrak{h}^*$ the set 
$R_{\lambda}=\{\alpha\in R\,:\, \lambda(H_{\alpha})\in\mathbb{Z}\}$ is a root system 
and the triangular decomposition of $\mathfrak{g}$ induces a uniquely defined 
basis, $B_{\lambda}$, of $R_{\lambda}$. Let $W_{\lambda}$ be the Weyl group of 
$R_{\lambda}$. We call $\lambda$ {\em relatively dominant} provided that $\lambda$ is 
a dominant  element in $\{w\cdot \lambda\,:\, w\in W_{\lambda}\}$ and 
{\em regular} provided that the stabilizer of $\lambda$ in $W_{\lambda}$ with 
respect to the dot action is trivial. 
\vspace{4mm}

{\em Throughout the paper we fix a relatively dominant and regular
$\lambda\in \mathfrak{h}^*$.}
\vspace{4mm}

For $w\in W_{\lambda}$ we denote by $\Delta(w)$ the Verma module with the highest 
weight $w\cdot\lambda$, and by $L(w)$ the unique simple quotient of $\Delta(w)$, 
see \cite[Chapter~7]{Di}. For $S\subset B_{\lambda}$ we denote by  $W^{S}_{\lambda}$ the 
subgroup of $W_{\lambda}$,  generated by $s_{\alpha}$, $\alpha\in S$. Denote by 
$w_{\lambda}^S$ the longest element in  $W^{S}_{\lambda}$
(in particular, $w_{\lambda}^{B_{\lambda}}$ is the longest element in $W_{\lambda}$). 
The main result of the present  paper is the following statement:

\begin{theorem}\label{tmain}
Let $S\subset B_{\lambda}$, $\alpha\in S$, and set 
$\mathbf{w}=s_{\alpha} w_{\lambda}^S w_{\lambda}^{B_{\lambda}}$. Then the canonical inclusion
\begin{displaymath}
U(\mathfrak{g})/\mathrm{Ann}(L(\mathbf{w}))\hookrightarrow\LL(L(\mathbf{w}),L(\mathbf{w}))
\end{displaymath}
is surjective.
\end{theorem}

The paper is organized as follows: in Section~\ref{s2} we collect all necessary
preliminaries on the category $\mathcal{O}$ and Arkhipov's twisting functors.
In Section~\ref{s3} we show how one can apply the twisting functors to obtain
the classical results related to Kostant's problem (in principal, if one takes into
account the relation between the twisting functors and Joseph's completion functors, 
obtained in \cite{KM}, our approach here is rather similar to the original approach,
however, here it is formulated in a shorter way). In Section~\ref{s4} we prove
Theorem~\ref{tmain} in the case $S=B_{\lambda}$. This is then
used in Section~\ref{s5} to prove Theorem~\ref{tmain} in the general case. In
Section~\ref{s6} we present an application of Theorem~\ref{tmain} and answer
Kostant's question for some simple $\alpha$-stratified modules.

\section{Preliminaries about the category $\mathcal{O}$}\label{s2}

Let $\mathcal{O}$ denote the BGG-category $\mathcal{O}$, associated with the
triangular decomposition of $\mathfrak{g}$, fixed above, see \cite{BGG}. Let 
${}^{\star}:\mathcal{O}\to \mathcal{O}$ be the classical duality on $\mathcal{O}$, 
that is a contravariant exact involutive equivalence, preserving the isomorphism 
classes of simple module, see \cite[Section~5]{Ir}. 
Let $\mathcal{O}_{\lambda}$ denote the indecomposable
block of $\mathcal{O}$, whose simple modules have the form $L(w)$, $w\in W_{\lambda}$.
Denote further by $P(w)$ the indecomposable projective cover of $L(w)$, see \cite{BGG}, 
and by $\theta_w$ the indecomposable projective functor on $\mathcal{O}_{\lambda}$,
uniquely determined by the property $\theta_w\Delta(e)\cong P(w)$, see \cite[I.3]{BG}.
Then $\{\theta_w\,:\,w\in W_{\lambda}\}$ are exactly the direct summands of the 
composition of $V\otimes {}_-$ followed by the projection from $\mathcal{O}$ 
to $\mathcal{O}_{\lambda}$, if we let $V$ run through all the finite-dimensional 
$\mathfrak{g}$-modules.

For $w\in W_{\lambda}$ set $\nabla(w)=\Delta(w)^{\star}$ and denote by 
$\mathcal{F}_{\lambda}(\Delta)$ the full subcategory of $\mathcal{O}_{\lambda}$, which 
consists of all modules, having a filtration, whose subquotients are isomorphic to 
Verma modules. Set $\mathcal{F}_{\lambda}(\nabla)=\mathcal{F}_{\lambda}(\Delta)^{\star}$.

Let $\mathcal{D}^b(\mathcal{O}_{\lambda})$ denote the bounded derived category 
of $\mathcal{O}_{\lambda}$. For a right or a left exact functor, $F$, on 
$\mathcal{O}_{\lambda}$ we denote by $\mathcal{L}F$ and $\mathcal{R}F$ the 
corresponding left and right derived functors respectively. For $i\geq 0$ we denote 
by $\mathcal{L}_iF$ and $\mathcal{R}^iF$ the corresponding $i$-th cohomology
functors. We denote by $[1]$ the shifting functor on 
$\mathcal{D}^b(\mathcal{O}_{\lambda})$ such that for every complex 
$\mathcal{X}^{\bullet}\in \mathcal{D}^b(\mathcal{O}_{\lambda})$
and for all $i\in\mathbb{Z}$ we have $\mathcal{X}[1]^{i}=\mathcal{X}^{i+1}$.
We consider $\mathcal{O}_{\lambda}$ as a subcategory of 
$\mathcal{D}^b(\mathcal{O}_{\lambda})$ via the classical embedding in degree zero. 

Via the equivalence from \cite{So} for $w\in W_{\lambda}$ we can define on 
$\mathcal{O}_{\lambda}$ {\em Arkhipov's twisting functor} $\mathtt{T}_w$, see 
\cite{Ar,AS,KM}, and denote by $\mathtt{G}_w$ its right adjoint (which is isomorphic, 
by \cite[Corollary~6]{KM}, to corresponding Joseph's completion functor from \cite{Jo1}, 
and to the functor $\star\mathtt{T}_w\star$, see \cite[Theorem~4.1]{AS}). In this 
paper we will use the following properties of $\mathtt{T}_w$ (the functor $\mathtt{G}_w$ 
has dual properties):

\begin{enumerate}[(I)]
\item\label{twist1} For every $w,x\in W_{\lambda}$ we have
$\mathtt{T}_w\theta_x\cong \theta_x\mathtt{T}_w$, see \cite[Theorem~3.2]{AS}.
\item\label{twist2} For every $w,x\in W_{\lambda}$ and $i>0$ we have 
$\mathcal{L}_i\mathtt{T}_w\Delta(x)=0$, see \cite[Theorem~2.2]{AS}.
\item\label{twist3} For every $w \in W_{\lambda}$ the functor $\mathcal{L}\mathtt{T}_w$
is an autoequivalence of $\mathcal{D}^b(\mathcal{O}_{\lambda})$ with the inverse 
functor $\mathcal{R}\mathtt{G}_{w^{-1}}$, see \cite[Corollary~4.2]{AS}.
\item\label{twist4} For every $w \in W_{\lambda}$ and every reduced decomposition,
$w=s_1\dots s_k$, we have $\mathtt{T}_w\cong \mathtt{T}_{s_1}\cdot \dots 
\cdot\mathtt{T}_{s_k}$, see \cite[Lemma~2.1]{AS} and \cite[Corollary~11]{KM}.
\item\label{twist5} For every $x \in W_{\lambda}$ and every simple reflection 
$s\in W_{\lambda}$ such that $sx>x$ we have $\mathtt{T}_{s}\Delta(x)\cong \Delta(sx)$,
see \cite[Lemma~6.2]{AL};
\item\label{twist6} For every $x \in W_{\lambda}$ and every simple reflection 
$s\in W_{\lambda}$ we have
\begin{displaymath}
\mathtt{T}_{s}\nabla(x)\cong 
\begin{cases}
\nabla(x), & x<sx,\\
\nabla(sx), & x>sx,
\end{cases}
\end{displaymath}
see \cite[Theorem~2.3]{AS}.
\item\label{twist7} For every $x \in W_{\lambda}$ and every simple reflection 
$s\in W_{\lambda}$ we have that $\mathtt{T}_{s}L(x)\neq 0$ if and only if $sx<x$,
see \cite[Section~6]{AS}.
\item\label{twist9} For every simple reflection $s\in W_{\lambda}$ and for every
$M\in\mathcal{O}_{\lambda}$ the module $\mathcal{L}_1\mathtt{T}_s(M)$ is the
largest $s$-finite submodule of $M$, see \cite[Theorem~1]{MS} or
\cite[Proposition~6]{Kh}.
\end{enumerate}

\section{The classical results}\label{s3}

We start with some preparation, during which we use the twisting functors to
obtain several classical results related to Kostant's problem. We base our 
approach on two classical statements. The first one, which can be found in 
\cite[6.8]{Ja}, is  a very abstract property of $\LL(M,N)$:

\begin{proposition}\label{prjan1}
Let $M,N$ be $\mathfrak{g}$-modules and $V$ be a finite-dimensional
$\mathfrak{g}$-module. Then there are canonical isomorphisms
\begin{equation}\label{eq241}
\mathrm{Hom}_{\mathfrak{g}}(V,\LL(M,N))\cong
\mathrm{Hom}_{\mathfrak{g}}(M\otimes V,N)\cong
\mathrm{Hom}_{\mathfrak{g}}(M,N\otimes V^*),
\end{equation}
where $\LL(M,N)$ is considered as a $\mathfrak{g}$-module
under the adjoint action.
\end{proposition}

The second statement is the classical positive answer to Kostant's problem
for projective Verma modules. In \cite[6.9]{Ja} it is shown that 

\begin{proposition}\label{prjan2}
For every submodule $M\subset \Delta(e)$ the canonical inclusion
\begin{displaymath}
U(\mathfrak{g})/\mathrm{Ann}(\Delta(e)/M)
\hookrightarrow
\LL(\Delta(e)/M,\Delta(e)/M)
\end{displaymath}
is surjective, in particular, the canonical inclusion
\begin{displaymath}
U(\mathfrak{g})/\mathrm{Ann}(\Delta(e))
\hookrightarrow
\LL(\Delta(e),\Delta(e))
\end{displaymath}
is surjective.
\end{proposition}

Using the twisting functors we obtain:

\begin{corollary}\label{cprjan2} (\cite[Corollary~6.4]{Jo2}, \cite[7.25]{Ja})
For every $w\in W_{\lambda}$  the canonical inclusion
\begin{displaymath}
U(\mathfrak{g})/\mathrm{Ann}(\Delta(w))
\hookrightarrow
\LL(\Delta(w),\Delta(w))
\end{displaymath}
is surjective.
\end{corollary}

\begin{proof}
We have the obvious map $\LL(\Delta(w),\Delta(w))\to 
\LL(\Delta(e),\Delta(e))$ induced by the inclusion
$\Delta(w)\subset \Delta(e)$. Since 
$\mathrm{Ann}(\Delta(w))=\mathrm{Ann}(\Delta(e))$ by
\cite[Theorem~8.4.4]{Di}, it is enough to show that for every simple
finite-dimensional $\mathfrak{g}$-module $V$ we have the equality
\begin{displaymath}
\dim\mathrm{Hom}_{\mathfrak{g}}(V,\LL(\Delta(w),\Delta(w)))=
\dim\mathrm{Hom}_{\mathfrak{g}}(V,\LL(\Delta(e),\Delta(e))).
\end{displaymath}
For this we compute
\begin{displaymath}
\begin{array}{lcl}
\mathrm{Hom}_{\mathfrak{g}}(V,\LL(\Delta(w),\Delta(w))) & = & 
\eqref{eq241}\\
\mathrm{Hom}_{\mathfrak{g}}(\Delta(w),\Delta(w)\otimes V^*) & = & 
\text{\eqref{twist4} and \eqref{twist5}}\\
\mathrm{Hom}_{\mathfrak{g}}(\mathrm{T}_w\Delta(e),\mathrm{T}_w(\Delta(e))\otimes V^*) 
& = &  \eqref{twist1}\\
\mathrm{Hom}_{\mathfrak{g}}(\mathrm{T}_w\Delta(e),\mathrm{T}_w(\Delta(e)\otimes V^*)) 
& = &  \\
\mathrm{Hom}_{\mathcal{D}^b(\mathcal{O}_{\lambda})}
(\mathrm{T}_w\Delta(e),\mathrm{T}_w(\Delta(e)\otimes V^*)) 
& = &  \eqref{twist2}\\
\mathrm{Hom}_{\mathcal{D}^b(\mathcal{O}_{\lambda})}(\mathcal{L}\mathrm{T}_w\Delta(e),
\mathcal{L}\mathrm{T}_w(\Delta(e)\otimes V^*)) & = &  \eqref{twist3}\\
\mathrm{Hom}_{\mathcal{D}^b(\mathcal{O}_{\lambda})}(\Delta(e),
\Delta(e)\otimes V^*) & = &  \\
\mathrm{Hom}_{\mathfrak{g}}(\Delta(e),
\Delta(e)\otimes V^*) & = & \eqref{eq241} \\
\mathrm{Hom}_{\mathfrak{g}}(V,\LL(\Delta(e),\Delta(e))). &  & 
\end{array}
\end{displaymath}
This completes the proof.
\end{proof}

\begin{proposition}\label{pr200}
Let $w\in W_{\lambda}$ and
\begin{equation}\label{eq444}
0\to X\to \Delta(w)\to Y\to 0
\end{equation}
be a short exact sequence such that for every finite-dimensional
$\mathfrak{g}$-module $V$ we have
\begin{equation}\label{eq445}
\mathrm{Ext}_{\mathcal{O}}^1(\Delta(w),X\otimes V)=0.
\end{equation}
Then the canonical inclusion
\begin{displaymath}
U(\mathfrak{g})/\mathrm{Ann}(Y)\hookrightarrow  \LL(Y,Y)
\end{displaymath}
is surjective.
\end{proposition}

\begin{proof}
Applying $\mathrm{Hom}_{\mathfrak{g}}(\Delta(w),{}_-\otimes V)$ to \eqref{eq444} and 
using \eqref{eq445} yields the short exact sequence
\begin{multline*}
0\to \mathrm{Hom}_{\mathfrak{g}}(\Delta(w),X\otimes V)\to
\mathrm{Hom}_{\mathfrak{g}}(\Delta(w),\Delta(w)\otimes V)\to\\ \to
\mathrm{Hom}_{\mathfrak{g}}(\Delta(w),Y\otimes V)\to 0,
\end{multline*}
which implies that $\LL(\Delta(w),\Delta(w))$ surjects onto
$\LL(\Delta(w),Y)$, where the vectorspace $\LL(Y,Y)$ is a subspace. 
Since $U(\mathfrak{g})$ surjects onto $\LL(\Delta(w),\Delta(w))$
by Corollary~\ref{cprjan2}, the statement follows.
\end{proof}

Now we can prove the classical result of Gabber and Joseph:

\begin{theorem}\label{c3}(\cite[Theorem~4.4]{GJ}, \cite[7.32]{Ja})
Let $S\subset B_{\lambda}$ and $\overline{\mathbf{w}}=w_{\lambda}^Sw_{\lambda}^{B_{\lambda}}$. 
Then the canonical inclusion  
\begin{displaymath}
U(\mathfrak{g})/\mathrm{Ann}(L(\overline{\mathbf{w}}))\hookrightarrow 
\LL(L(\overline{\mathbf{w}}),L(\overline{\mathbf{w}}))
\end{displaymath}
 is surjective. 
\end{theorem}

\begin{proof}
Let $V$ be a finite-dimensional $\mathfrak{g}$-module. Consider the short exact sequence
\begin{equation}\label{eq333}
0\to K(\overline{\mathbf{w}})\to \Delta(\overline{\mathbf{w}})\to L(\overline{\mathbf{w}})\to 0,
\end{equation}
where $K(\overline{\mathbf{w}})$ is just the kernel of the canonical projection from 
$\Delta(\overline{\mathbf{w}})$ to $L(\overline{\mathbf{w}})$.
Then we have
\begin{displaymath}
\begin{array}{lll}
\mathrm{Ext}_{\mathcal{O}}^1(\Delta(\overline{\mathbf{w}}),K(\overline{\mathbf{w}})\otimes V) & = & \\
\mathrm{Hom}_{\mathcal{D}^b(\mathcal{O})}
(\Delta(\overline{\mathbf{w}}),K(\overline{\mathbf{w}})\otimes V[1]) & = & 
\text{(duals of \eqref{twist4} and \eqref{twist5})}\\
\mathrm{Hom}_{\mathcal{D}^b(\mathcal{O})}
(\mathtt{G}_{w_{\lambda}^S}\Delta(w_{\lambda}^{B_{\lambda}}),K(\overline{\mathbf{w}})\otimes V[1]) & = & 
\text{(${}_-\otimes V$ is exact}\\
 & & \text{and preserves projectives)}\\
\mathrm{Hom}_{\mathcal{D}^b(\mathcal{O})}
(V^{*}\otimes \mathtt{G}_{w_{\lambda}^S}\Delta(w_{\lambda}^{B_{\lambda}}),K(\overline{\mathbf{w}})[1]) & = & 
\text{(dual of \eqref{twist1})}\\
\mathrm{Hom}_{\mathcal{D}^b(\mathcal{O})}
(\mathtt{G}_{w_{\lambda}^S}(V^{*}\otimes \Delta(w_{\lambda}^{B_{\lambda}})),K(\overline{\mathbf{w}})[1]) & = & 
\text{($\Delta(w_{\lambda}^{B_{\lambda}})=\nabla(w_{\lambda}^{B_{\lambda}})$)}\\
\mathrm{Hom}_{\mathcal{D}^b(\mathcal{O})}
(\mathtt{G}_{w_{\lambda}^S}(V^{*}\otimes \nabla(w_{\lambda}^{B_{\lambda}})),K(\overline{\mathbf{w}})[1]) & = & 
\text{(dual of \eqref{twist2})}\\
\mathrm{Hom}_{\mathcal{D}^b(\mathcal{O})}
(\mathcal{R}\mathtt{G}_{w_{\lambda}^S}(V^{*}\otimes 
\nabla(w_{\lambda}^{B_{\lambda}})),K(\overline{\mathbf{w}})[1]) & = & 
\text{\eqref{twist3}}\\
\mathrm{Hom}_{\mathcal{D}^b(\mathcal{O})}
(V^{*}\otimes \nabla(w_{\lambda}^{B_{\lambda}}),
\mathcal{L}\mathtt{T}_{w_{\lambda}^S}K(\overline{\mathbf{w}})[1]). 
&&
\end{array}
\end{displaymath}
Let us calculate $\mathcal{L}\mathtt{T}_{w_{\lambda}^S}K(\overline{\mathbf{w}})$. 
Because of our choice of $w_{\lambda}^S$ we can use Proposition~\ref{preqv},
which will be proved in Section~\ref{s5} (alternatively one can use \cite[Section~2]{GJ}),
and \cite{BGG2} to get that the module $K(\overline{\mathbf{w}})$ admits a BGG-type 
resolution, which has the following form:
\begin{displaymath}
0\to X^{k}\to X^{k-1}\to \dots\to X^1\to X^0\to K(\overline{\mathbf{w}})\to 0,
\end{displaymath}
where every $X^i$ is a direct sum of some $\Delta(xw_{\lambda}^{B_{\lambda}})$ with
$x\in W_{\lambda}^S$, $x\neq w_{\lambda}^S$. Let  $\mathcal{X}^{\bullet}$ denote the 
corresponding complex in $\mathcal{D}^b(\mathcal{O}_{\lambda})$. Using
\eqref{twist2} we have
\begin{displaymath}
\mathcal{L}\mathtt{T}_{w^S_{\lambda}}K(\overline{\mathbf{w}})=
\mathcal{L}\mathtt{T}_{w^S_{\lambda}}\mathcal{X}^{\bullet}=
\mathtt{T}_{w^S_{\lambda}}\mathcal{X}^{\bullet}.
\end{displaymath}
Now for every $x\in W_{\lambda}^S$, $x\neq w_{\lambda}^S$, 
let $y\in W_{\lambda}^S$ be such that $yx^{-1}=w_{\lambda}^S$. Then, using 
\eqref{twist4}, \eqref{twist5}, and \eqref{twist6},  we  have
\begin{displaymath}
\mathtt{T}_{w_{\lambda}^S}\Delta(xw_{\lambda}^{B_{\lambda}})=
\mathtt{T}_{y}\mathtt{T}_{x^{-1}}\Delta(xw_{\lambda}^{B_{\lambda}})=
\mathtt{T}_{y}\Delta(w_{\lambda}^{B_{\lambda}})=
\mathtt{T}_{y}\nabla(w_{\lambda}^{B_{\lambda}})=
\nabla(yw_{\lambda}^{B_{\lambda}}).
\end{displaymath}
This implies that $\mathtt{T}_{w_{\lambda}^S}\mathcal{X}^{\bullet}$ is a complex of dual
Verma modules in $\mathcal{O}_{\lambda}$. At the same time the module 
$\nabla(w_{\lambda}^{B_{\lambda}})\otimes V^*\cong
\Delta(w_{\lambda}^{B_{\lambda}})\otimes V^*$  is a tilting module in 
$\mathcal{O}_{\lambda}$. Hence, by \cite[Chap.~III,~Lemma~2.1]{Ha}, the space 
\begin{displaymath}
\mathrm{Hom}_{\mathcal{D}^b(\mathcal{O})}
(\nabla(w_{\lambda}^{B_{\lambda}})\otimes V^*,
\mathcal{L}\mathtt{T}_{w_{\lambda}^S}K(\overline{\mathbf{w}})[1])
\end{displaymath}
can be computed already in the homotopy category, where it is obviously zero, since the 
only non-zero component of the first complex is in degree zero and the above computation
shows that the zero component of the second complex is zero. Hence we obtain
\begin{equation}\label{eq334}
\mathrm{Ext}_{\mathcal{O}}^1(\Delta(\overline{\mathbf{w}}),K(\overline{\mathbf{w}})\otimes V)=0.
\end{equation}
The statement of our theorem now follows by applying Proposition~\ref{pr200}
to the short exact sequence \eqref{eq333}.
\end{proof}

\section{Proof of Theorem~\ref{tmain}: the case $S=B_{\lambda}$}\label{s4}

In this section we prove Theorem~\ref{tmain} in the case $S=B_{\lambda}$. 
Throughout the section we fix $\alpha\in B_{\lambda}$ and set $s=s_{\alpha}$.

\begin{proposition}\label{pr10}
The canonical inclusion $U(\mathfrak{g})/\mathrm{Ann}(L(s)) \hookrightarrow
\LL(L(s), L(s))$ is surjective.
\end{proposition}

To prove this statement we will need several lemmas.

\begin{lemma}\label{l4.1}
Let $w\in W_{\lambda}$ be such that 
$\mathrm{Hom}_{\mathfrak{g}}(L(s), \theta_wL(s))\neq 0$.
Then $w=s$ or $w=e$.
\end{lemma}

\begin{proof}
Assume that $w\neq e,s$. Let
\begin{equation}\label{eq340}
0\to L(s)\to X\to L(e)\to 0
\end{equation}
be a non-split short exact sequence, which exists because of the Kazhdan-Lusztig
theorem (see for example \cite[Theorem~1]{Ko}). 
Then $\theta_w L(e)=0$ since $w\neq e$ and hence
$\theta_w X=\theta_w L(s)$. However, since \eqref{eq340} is non-split,
$X$ is a homomorphic image of $\Delta(e)$, and hence 
$\theta_w X$ is a homomorphic image of $\theta_w\Delta(e)=P(w)$. In particular,
$\theta_w X$ is either zero or has simple top $L(w)$. On the other hand 
$\theta_w L(s)$ is self-dual and thus $\theta_w X=\theta_w L(s)$ is either zero
or has simple socle $L(w)$. In each of these two cases we have the equality
$\mathrm{Hom}_{\mathfrak{g}}(L(s), \theta_w L(s))= 0$
since $w\neq s$. This completes the proof.
\end{proof}

The above result naturally motivates the following question:
\vspace{3mm}

{\em {\bf Question.} Let $S\subset B_{\lambda}$ and $w\in W_{\lambda}$ be such that 
the vector space $\mathrm{Hom}_{\mathfrak{g}}(L(w_{\lambda}^{S}), \theta_w
L(w_{\lambda}^{S}))$ is non-zero. Does this imply that 
$w\in W_{\lambda}^S$?}
\vspace{3mm}

Recall that a $\mathfrak{g}$-module, $M$, is called {\em $s$-finite} provided that
it is locally finite over the $\mathfrak{sl}_2$-subalgebra of $\mathfrak{g}$,
which corresponds to $s$. The module $L(w)$ is $s$-finite if and only if
$w$ is the minimal coset representative of some coset from 
${\{e,s\}}\backslash W_{\lambda}$, that is if and only if $sw>w$.

Define $F(s)$ as the minimal submodule of the radical $\mathrm{Rad}(\Delta(s))$
of $\Delta(s)$ such that the quotient $\mathrm{Rad}(\Delta(s))/F(s)$ is
$s$-finite and consider the short exact sequence
\begin{equation}\label{eq432}
0\to F(s)\to \Delta(s)\to N(s)\to 0,
\end{equation}
where $N(s)$ is the cokernel. Our next step is to prove the following:

\begin{lemma}\label{l4.2}
The canonical inclusion $U(\mathfrak{g})/\mathrm{Ann}(N(s)) \hookrightarrow
\LL(N(s), N(s))$ is surjective.
\end{lemma}

\begin{proof}
For every $w\in W_{\lambda}$ we have
\begin{displaymath}
\begin{array}{lcl}
\mathrm{Ext}_{\mathcal{O}}^1(\Delta(s),\theta_w F(s)) & = & \\
\mathrm{Hom}_{\mathcal{D}^b(\mathcal{O})}(\Delta(s),\theta_w F(s)[1]) & = & 
\text{\eqref{twist5}}\\
\mathrm{Hom}_{\mathcal{D}^b(\mathcal{O})}(\mathtt{T}_s\Delta(e),\theta_w F(s)[1]) & = & 
\text{(properties of $\theta_{w}$)}\\
\mathrm{Hom}_{\mathcal{D}^b(\mathcal{O})}(\theta_{w^{-1}}\mathtt{T}_s\Delta(e), F(s)[1]) 
& = &  \text{\eqref{twist1}}\\
\mathrm{Hom}_{\mathcal{D}^b(\mathcal{O})}(\mathtt{T}_s\theta_{w^{-1}}\Delta(e), F(s)[1]) 
& = &  \\
\mathrm{Hom}_{\mathcal{D}^b(\mathcal{O})}(\mathtt{T}_s P(w^{-1}), F(s)[1]) 
& = &  \text{\eqref{twist2}}\\
\mathrm{Hom}_{\mathcal{D}^b(\mathcal{O})}(\mathcal{L}\mathtt{T}_s P(w^{-1}), F(s)[1]) 
& = &  \text{\eqref{twist3}}\\
\mathrm{Hom}_{\mathcal{D}^b(\mathcal{O})}(P(w^{-1}), \mathcal{R}\mathtt{G}_sF(s)[1]). 
&=  & \text{($P(w^{-1})$ is projective)}\\
\mathrm{Hom}_{\mathfrak{g}}(P(w^{-1}), \mathcal{R}^1\mathtt{G}_sF(s)) & =&
\text{(dual of \eqref{twist9})}\\
0. & & 
\end{array}
\end{displaymath}
The statement now follows from Corollary~\ref{cprjan2} and Proposition~\ref{pr200}.
\end{proof}

Consider now the short exact sequence
\begin{equation}\label{eq456}
0\to X(s)\to N(s)\overset{p'}{\longrightarrow} L(s)\to 0.
\end{equation}

\begin{lemma}\label{l4.3}
For every finite-dimensional $\mathfrak{g}$-module $V$ the sequence
\eqref{eq456} induces the following isomorphism:
\begin{displaymath}
\mathrm{Hom}_{\mathfrak{g}}(N(s),N(s)\otimes V)\cong
\mathrm{Hom}_{\mathfrak{g}}(L(s),L(s)\otimes V).
\end{displaymath}
\end{lemma}

\begin{proof}
Let $0\neq f\in \mathrm{Hom}_{\mathfrak{g}}(N(s),N(s)\otimes V)$. Since
$N(s)$ has simple top, $f$ can not annihilate it. Consider the map
$(p'\otimes\mathrm{id})\circ f\in \mathrm{Hom}_{\mathfrak{g}}(N(s),L(s)\otimes V)$. 
Since the kernel of the projection 
$N(s)\otimes V\overset{p'\otimes \mathrm{id}}{\tto}L(s)\otimes V$
is $s$-finite and the top of $N(s)$ is not, we have $(p'\otimes\mathrm{id})\circ f\neq 0$.
On the other hand, since the socle of $L(s)\otimes V$ consists exclusively of
$s$-infinite modules, the map  $(p'\otimes\mathrm{id})\circ f$ must annihilate $X(s)$ and 
hence it factors through $L(s)$. This implies that \eqref{eq456} induces the following 
inclusion:
\begin{displaymath}
\mathrm{Hom}_{\mathfrak{g}}(N(s),N(s)\otimes V)\hookrightarrow
\mathrm{Hom}_{\mathfrak{g}}(L(s),L(s)\otimes V).
\end{displaymath}
To complete the proof we now have to compare the dimensions and thus it
is enough to show that for every $w\in W_{\lambda}$ we have
\begin{displaymath}
\dim\mathrm{Hom}_{\mathfrak{g}}(N(s),\theta_w N(s))=
\dim\mathrm{Hom}_{\mathfrak{g}}(L(s),\theta_w L(s)).
\end{displaymath}
This is obvious for $w=e$ since both spaces are one-dimensional in this case. 
For $w=s$ we have non-zero adjunction morphisms in
both spaces, moreover, the module $\theta_s L(s)$ has simple socle. This implies
\begin{displaymath}
\dim\mathrm{Hom}_{\mathfrak{g}}(N(s),\theta_s N(s))=
\dim\mathrm{Hom}_{\mathfrak{g}}(L(s),\theta_s L(s))=1.
\end{displaymath}
For $w\neq s,e$ Lemma~\ref{l4.1} implies 
$\dim\mathrm{Hom}_{\mathfrak{g}}(L(s),\theta_w L(s))=0$.
The statement follows.
\end{proof}

Now we are ready to prove Proposition~\ref{pr10}.

\begin{proof}[Proof of Proposition~\ref{pr10}.]
Since $X(s)$ is $s$-finite and $L(s)$ is simple and $s$-in\-fi\-ni\-te, we have
$\LL(X(s),L(s))=0$ by Proposition~\ref{prjan1}, which implies that \eqref{eq456}
induces the isomorphism $\LL(L(s),L(s))\cong \LL(N(s),L(s))$.

Since $X(s)$ is $s$-finite and the top of $N(s)$ is simple and $s$-infinite, we have
$\LL(N(s),X(s))=0$ by Proposition~\ref{prjan1}, which implies that \eqref{eq456}
induces the inclusion $\LL(N(s),N(s))\hookrightarrow  \LL(N(s),L(s))$.
However, Lemma~\ref{l4.3} and Proposition~\ref{prjan1} show that this inclusion
is in fact an isomorphism. Since $U(\mathfrak{g})$ surjects onto
$\LL(N(s),N(s))$ by Corollary~\ref{cprjan2}, it follows that \eqref{eq456}
induces a surjection of $U(\mathfrak{g})$ onto $\LL(L(s),L(s))$. This 
completes the proof.
\end{proof}

\section{Proof of Theorem~\ref{tmain}: the general case}\label{s5}

In this Section we prove Theorem~\ref{tmain} in the general case. Our approach is similar
to the one we use in Section~\ref{s4}, however, it requires more delicate arguments in
several places, moreover, in some places we will use the reduction to the case,
considered in Section~\ref{s4}. Set $s=s_{\alpha}$ and recall the notation
$\overline{\mathbf{w}}=w_{\lambda}^Sw_{\lambda}^{B_{\lambda}}$.

Using the equivalence from \cite{So} we can assume that $\lambda$ is integral.
Let $\mathfrak{a}=\mathfrak{a}(S)$ denote the
semi-simple Lie subalgebra of $\mathfrak{g}$, generated by $X_{\pm\alpha}$, $\alpha\in S$.
If $M$ is a weight $\mathfrak{g}$-module with the weight-space decomposition
$M=\oplus_{\mu\in\mathfrak{h}^*}M_{\mu}$, and $\nu\in \mathfrak{h}^*$, then the subspace
\begin{displaymath}
M_{\mathfrak{a}}^{\nu}=\oplus_{\mu\in \nu+\mathbb{Z}S}M_{\mu}
\end{displaymath}
is stable under the action of $\mathfrak{a}$ and hence is an $\mathfrak{a}$-submodule of $M$.
This induces the functor, which we will denote by  $\mathtt{R}^{\nu}$, from the category of 
all weight $\mathfrak{g}$-modules to the category of  all weight $\mathfrak{a}$-modules,
which sends $M$ to $M_{\mathfrak{a}}^{\nu}$. Let $\mathcal{O}^{\mathfrak{a}}$ denote the category
$\mathcal{O}$ for the algebra $\mathfrak{a}$. From the PBW theorem it follows that 
for every $w\in W_{\lambda}$ and every $\nu\in \mathfrak{h}^*$ the module
$\mathtt{R}^{\nu}\Delta(w)$ has a finite Verma flag as an $\mathfrak{a}$-module,
in particular, $\mathtt{R}^{\nu}\Delta(w)\in \mathcal{O}^{\mathfrak{a}}$.
From this one easily deduces that $\mathtt{R}^{\nu}$ maps $\mathcal{O}$ to
$\mathcal{O}^{\mathfrak{a}}$.

Let $\mathfrak{h}^{\perp}$ be the orthogonal complement to 
$\mathfrak{a}\cap \mathfrak{h}$ in $\mathfrak{h}$ with respect to the Killing
form. Let $\xi$ be the restriction of $\overline{\mathbf{w}}\cdot \lambda$
to $\mathfrak{h}^{\perp}$. Define the {\em parabolic induction functor}
$\mathrm{Ind}_{\mathfrak{a}}^{\mathfrak{g}}$ in the following way: for
$M\in \mathcal{O}^{\mathfrak{a}}$ let $\mathfrak{h}^{\perp}$ act on 
$M$ via $\xi$,  and let $X_{\alpha}M=0$ for all positive roots $\alpha\in R$
such that $X_{\alpha}\not\in \mathfrak{a}$. In this way we can regard
$M$ as a module over the parabolic subalgebra
$\mathfrak{p}=\mathfrak{a}+\mathfrak{h}+\mathfrak{n}_+$ of $\mathfrak{g}$. 
We set
\begin{displaymath}
\mathrm{Ind}_{\mathfrak{a}}^{\mathfrak{g}}(M)=U(\mathfrak{g})
\otimes_{U(\mathfrak{p})}M,
\end{displaymath}
which obviously defines a functor from $\mathcal{O}^{\mathfrak{a}}$ to
$\mathcal{O}$. From the PBW theorem it follows that this functor sends
Verma modules to Verma modules. Let $\zeta$ be the restriction of
$\overline{\mathbf{w}}\cdot \lambda$ to 
$\mathfrak{a}\cap \mathfrak{h}$. Note that $\zeta$ is regular and dominant
for $\mathfrak{a}$. 

Finally, denote by $\mathcal{C}$ the full subcategory of $\mathcal{O}^{\lambda}$,
which consists of all modules $M$, whose all composition factors have the form
$L(y)$, $y\in W_{\lambda}^Sw_{\lambda}^{B_{\lambda}}$. 

\begin{proposition}\label{preqv}
$\mathrm{Ind}_{\mathfrak{a}}^{\mathfrak{g}}$ and
$\mathtt{R}^{\overline{\mathbf{w}}\cdot \lambda}$ induce mutually inverse
equivalences between $\mathcal{O}^{\mathfrak{a}}_{\zeta}$ and
$\mathcal{C}$.
\end{proposition}

\begin{proof}
The classical adjunction between the restriction and induction implies
that $(\mathrm{Ind}_{\mathfrak{a}}^{\mathfrak{g}},
\mathtt{R}^{\overline{\mathbf{w}}\cdot \lambda})$ is an adjoint pair
of functors, which gives us the natural maps
\begin{displaymath}
\mathrm{Ind}_{\mathfrak{a}}^{\mathfrak{g}}\,
\mathtt{R}^{\overline{\mathbf{w}}\cdot \lambda}\to
\mathrm{Id}_{\mathcal{C}},\quad\text{and}\quad
\mathrm{Id}_{\mathcal{O}^{\mathfrak{a}}_{\zeta}}\to
\mathtt{R}^{\overline{\mathbf{w}}\cdot \lambda}\,
\mathrm{Ind}_{\mathfrak{a}}^{\mathfrak{g}}.
\end{displaymath}
These maps are obviously isomorphisms on Verma modules, and then by induction
one shows that they are isomorphisms on simple modules. The statement follows.
\end{proof}

As in Section~\ref{s4}
we define $F(\mathbf{w})$ as the minimal submodule of the radical 
$\mathrm{Rad}(\Delta(\mathbf{w}))$  of $\Delta(\mathbf{w})$ such that the 
quotient $\mathrm{Rad}(\Delta(\mathbf{w}))/F(\mathbf{w})$ is
$s$-finite and consider the short exact sequence
\begin{equation}\label{eq43201}
0\to F(\mathbf{w})\to \Delta(\mathbf{w})\to N(\mathbf{w})\to 0,
\end{equation}
where $N(\mathbf{w})$ is the cokernel. 

\begin{proposition}\label{l5.1}
The canonical inclusion 
\begin{displaymath}
U(\mathfrak{g})/\mathrm{Ann}(N(\mathbf{w})) \hookrightarrow 
\LL(N(\mathbf{w}), N(\mathbf{w}))
\end{displaymath}
is surjective.
\end{proposition}

\begin{proof}
Let $w\in W_{\lambda}$. Using the same arguments as in the proof of
Lemma~\ref{l4.2} we obtain
\begin{displaymath}
\mathrm{Ext}_{\mathcal{O}}^1(\Delta(\mathbf{w}),\theta_w F(\mathbf{w}))=
\mathrm{Ext}_{\mathcal{O}}^1(\theta_{w^{-1}}
\Delta(\overline{\mathbf{w}}),\mathtt{G}_sF(\mathbf{w})).
\end{displaymath}
Let us prove that the last space is zero. For this we will
need the following statement:

\begin{lemma}\label{l4.5}
All simple subquotients of $\mathtt{G}_sF(\mathbf{w})$ are of the form $L(x)$,
$x\in W_{\lambda}^Sw_{\lambda}^{B_{\lambda}}$.
\end{lemma}

\begin{proof}
Using the left exactness of $\mathtt{G}_s$ it is enough to prove that 
for every $y\in W_{\lambda}^Sw_{\lambda}^{B_{\lambda}}$ all simple subquotients 
of $\mathtt{G}_s L(y)$ 
are of the form $L(x)$, $x\in W_{\lambda}^Sw_{\lambda}^{B_{\lambda}}$. 
By \eqref{twist7}, we can even assume $sy<y$.
Applying $\mathtt{G}_s$ to the short exact sequence
\begin{displaymath}
0\to K(y)\to \Delta(y)\to L(y)\to 0
\end{displaymath}
and using the dual of \eqref{twist6} we obtain the following exact sequence:
\begin{displaymath}
0\to \mathtt{G}_s K(y)\to \mathtt{G}_s\Delta(y)\left(\cong\Delta(sy)\right)
\to \mathtt{G}_sL(y)\to \mathcal{L}_1\mathtt{G}_s K(y).
\end{displaymath}
Obviously all simple subquotients of $\Delta(sy)$ have the necessary form and
from the dual of \eqref{twist9} it follows that 
all simple subquotients of $\mathcal{L}_1\mathtt{G}_s K(y)$
have the necessary form as well. The statement follows.
\end{proof}

We have $\theta_{w^{-1}}\Delta(\overline{\mathbf{w}})\in
\mathcal{F}_{\lambda}(\Delta)$. Let
$Q_1=\oplus_{x\in W_{\lambda}\setminus W_{\lambda}^Sw_{\lambda}^{B_{\lambda}}}P(x)$
and consider the trace $Z$ (i.e. the sum of the images of all homomorphism) 
of $Q_1$ in $\theta_{w^{-1}}\Delta(\overline{\mathbf{w}})$. 
Since the extensions between the Verma modules are directed, it follows by induction 
that all modules in the short exact sequence
\begin{displaymath}
0\to Z\to \theta_{w^{-1}}\Delta(\overline{\mathbf{w}})
\to\mathrm{Coker}\to 0
\end{displaymath}
have Verma flags, moreover, the Verma modules, occurring as subquotients
of $Z$, have the form $\Delta(x)$, $x\in W_{\lambda}\setminus 
W_{\lambda}^Sw_{\lambda}^{B_{\lambda}}$,
and the Verma modules, occurring as subquotients of 
$\mathrm{Coker}$, have the form $\Delta(y)$, $y\in 
W_{\lambda}^Sw_{\lambda}^{B_{\lambda}}$.
Since for every $x\in W_{\lambda}\setminus W_{\lambda}^Sw_{\lambda}^{B_{\lambda}}$ 
and $y\in W_{\lambda}^Sw_{\lambda}^{B_{\lambda}}$
we have $x\not\leq y$ with respect to the Bruhat order, for all such $x$ and $y$
we obtain 
\begin{displaymath}
\mathrm{Ext}_{\mathcal{O}}^1(\Delta(x),L(y))=0,
\end{displaymath}
which, because of Lemma~\ref{l4.5}, yields
\begin{displaymath}
\mathrm{Ext}_{\mathcal{O}}^1(Z,\mathtt{G}_sF(\mathbf{w}))=0.
\end{displaymath}

Now let us consider the module $\mathrm{Coker}\in\mathcal{C}$. We claim that 
$\mathtt{R}^{\overline{\mathbf{w}}\cdot\lambda}\mathrm{Coker}$ is a 
projective module in the category $\mathcal{O}^{\mathfrak{a}}_{\zeta}$. 
Indeed, the module $\Delta(\overline{\mathbf{w}})$ is
obtained by the parabolic induction from some projective Verma 
$\mathfrak{a}$-module. Since the adjoint action of $\mathfrak{a}$ on  
$U(\mathfrak{g})$ is locally finite, it follows that 
$\mathtt{R}^{\nu}(\Delta(\overline{\mathbf{w}}))$
is projective in $\mathcal{O}^{\mathfrak{a}}_{\zeta}$ for every
$\nu\in \mathfrak{h}^*$. Further, for every finite-dimensional
$\mathfrak{g}$-module $V$ we have 
\begin{displaymath}
\mathtt{R}^{\overline{\mathbf{w}}\cdot\lambda}
(V\otimes \Delta(\overline{\mathbf{w}}))=
\oplus_{(\nu_1,\nu_2)}
\mathtt{R}^{\nu_1}V\otimes
\mathtt{R}^{\nu_2}\Delta(\overline{\mathbf{w}}),
\end{displaymath}
where the sum is taken over all pairs $(\nu_1,\nu_2)\in\mathfrak{h}^*\times 
\mathfrak{h}^*$ with different $\mathfrak{h}^{\perp}$-restrictions of $\nu_1$ 
such that $\nu_1+\nu_2=\overline{\mathbf{w}}\cdot\lambda$. In particular,
$\mathtt{R}^{\overline{\mathbf{w}}\cdot\lambda}
(V\otimes \Delta(\overline{\mathbf{w}}))$ is projective
in $\mathcal{O}^{\mathfrak{a}}_{\zeta}$. The inductive
construction of the Verma flag in \cite{BGG} implies that 
\begin{displaymath}
\mathtt{R}^{\overline{\mathbf{w}}\cdot\lambda}
(\mathrm{Coker})= \mathtt{R}^{0}V\otimes
\mathtt{R}^{\overline{\mathbf{w}}\cdot\lambda}\Delta(\overline{\mathbf{w}}),
\end{displaymath}
which is also projective in $\mathcal{O}^{\mathfrak{a}}_{\zeta}$.
In particular, the first extension between
$\mathtt{R}^{\overline{\mathbf{w}}\cdot\lambda} (\mathrm{Coker})$
and all simple $\mathfrak{a}$-modules in 
$\mathcal{O}^{\mathfrak{a}}_{\zeta}$ vanishes and hence
from Proposition~\ref{preqv} we derive
\begin{displaymath}
\mathrm{Ext}_{\mathcal{O}}^1(\mathrm{Coker},L(y))=0
\end{displaymath}
for all $y\in W_{\lambda}^Sw_{\lambda}^{B_{\lambda}}$. Therefore, 
using Lemma~\ref{l4.5} we get
\begin{displaymath}
\mathrm{Ext}_{\mathcal{O}}^1(\mathrm{Coker},\mathtt{G}_sF(\mathbf{w}))=0.
\end{displaymath}

Thus 
\begin{displaymath}
\mathrm{Ext}_{\mathcal{O}}^1(\theta_{w^{-1}}
\Delta(\overline{\mathbf{w}}),\mathtt{G}_sF(\mathbf{w}))=0
\end{displaymath}
and the statement of the proposition follows from 
Corollary~\ref{cprjan2} and Proposition~\ref{pr200}.
\end{proof}

Consider now the short exact sequence
\begin{equation}\label{eq45601}
0\to X(\mathbf{w})\to N(\mathbf{w})\overset{\hat{p}}{\longrightarrow} L(\mathbf{w})\to 0.
\end{equation}

\begin{lemma}\label{l5.5}
For every finite-dimensional $\mathfrak{g}$-module $V$ 
the sequence \eqref{eq45601} induces the isomorphism
\begin{displaymath}
\mathrm{Hom}_{\mathfrak{g}}(N(\mathbf{w}),N(\mathbf{w})\otimes V)\cong
\mathrm{Hom}_{\mathfrak{g}}(L(\mathbf{w}),L(\mathbf{w})\otimes V).
\end{displaymath}
\end{lemma}

\begin{proof}
The same arguments as in Lemma~\ref{l4.3} show that \eqref{eq45601}
induces the inclusion
\begin{displaymath}
\mathrm{Hom}_{\mathfrak{g}}(N(\mathbf{w}),N(\mathbf{w})\otimes V)\hookrightarrow
\mathrm{Hom}_{\mathfrak{g}}(L(\mathbf{w}),L(\mathbf{w})\otimes V).
\end{displaymath}
Let $f\in \mathrm{Hom}_{\mathfrak{g}}(L(\mathbf{w}),L(\mathbf{w})\otimes V)$. We would
like to  lift $f$ to an element in the space
$\mathrm{Hom}_{\mathfrak{g}}(N(\mathbf{w}),N(\mathbf{w})\otimes V)$. For this
we consider the auxiliary module $\mathtt{G}_sL(\mathbf{w})$. 
 
\begin{lemma}\label{l5.6}
$\mathrm{Ext}_{\mathcal{O}}^1(\mathtt{G}_sL(\mathbf{w}),L(x))=0$ for 
each $s$-finite $L(x)$.
\end{lemma}

\begin{proof}
This follows from \cite[Theorem~6.3(3)]{AS} and the Kazhdan-Lusztig 
theorem (observe that in \cite[Theorem~6.3(3)]{AS} one has to assume that
$L'$ is $s$-finite).
\end{proof}

First of all we claim that $\mathtt{G}_sL(\mathbf{w})$ surjects onto $N(\mathbf{w})$. Indeed,
applying $\mathrm{Hom}_{\mathfrak{g}}(\mathtt{G}_sL(\mathbf{w}),{}_-)$ to
\eqref{eq45601} and using Lemma~\ref{l5.6} we obtain the surjection
\begin{displaymath}
\mathrm{Hom}_{\mathfrak{g}}(\mathtt{G}_sL(\mathbf{w}),N(\mathbf{w}))\tto
\mathrm{Hom}_{\mathfrak{g}}(\mathtt{G}_sL(\mathbf{w}),L(\mathbf{w})).
\end{displaymath}
Using this surjection we can lift the canonical projection
$\tilde{p}:\mathtt{G}_sL(\mathbf{w})\tto L(\mathbf{w})$ to obtain the short exact sequence
\begin{equation}\label{eq777}
0\to\mathrm{Ker}\to \mathtt{G}_sL(\mathbf{w})\overset{q}{\longrightarrow} N(\mathbf{w})\to 0.
\end{equation}

Applying $\mathrm{Hom}_{\mathfrak{g}}(\mathtt{G}_sL(\mathbf{w}),{}_-\otimes V)$ to
\eqref{eq45601} and using Lemma~\ref{l5.6} we obtain the surjection
\begin{displaymath}
\mathrm{Hom}_{\mathfrak{g}}(\mathtt{G}_sL(\mathbf{w}),N(\mathbf{w})\otimes V)\tto
\mathrm{Hom}_{\mathfrak{g}}(\mathtt{G}_sL(\mathbf{w}),L(\mathbf{w})\otimes V).
\end{displaymath}
In particular, we can lift the map 
$f\circ \hat{p}\circ q\in 
\mathrm{Hom}_{\mathfrak{g}}(\mathtt{G}_sL(\mathbf{w}),L(\mathbf{w})\otimes V)$
to some map $\overline{f}\in 
\mathrm{Hom}_{\mathfrak{g}}(\mathtt{G}_sL(\mathbf{w}),N(\mathbf{w})\otimes V)$.

Now recall that $\mathtt{G}_sL(\mathbf{w})\in\mathcal{C}$ by Lemma~\ref{l4.5}. Applying
$\mathtt{R}^{\overline{\mathbf{w}}\cdot\lambda}$ and using
Proposition~\ref{preqv} and Lemma~\ref{l4.3} we obtain that 
$\mathtt{R}^{\overline{\mathbf{w}}\cdot\lambda}(\overline{f})$ annihilates the module
$\mathtt{R}^{\overline{\mathbf{w}}\cdot\lambda}(\mathrm{Ker})$, which implies that 
$\overline{f}$ annihilates $\mathrm{Ker}$ by Proposition~\ref{preqv}. In particular,
$\overline{f}$ factors through $N(\mathbf{w})$. Since all the modules
$L(\mathbf{w})$, $N(\mathbf{w})$, and $\mathtt{G}_sL(\mathbf{w})$, have same simple top, it follows that
$f\neq 0$ if and only if $\overline{f}\neq 0$. This gives us the injection
\begin{displaymath}
\mathrm{Hom}_{\mathfrak{g}}(L(\mathbf{w}),L(\mathbf{w})\otimes V)\hookrightarrow
\mathrm{Hom}_{\mathfrak{g}}(N(\mathbf{w}),N(\mathbf{w})\otimes V),
\end{displaymath}
and the statement follows.
\end{proof}

Now we have the same amount of information as at the end of
Section~\ref{s4} and hence the proof of Theorem~\ref{tmain} can
be easily completed in the same way as the proof of Proposition~\ref{pr10}.

\begin{proof}[Proof of Theorem~\ref{tmain}.]
Mutatis mutandis the proof of Proposition~\ref{pr10}.
\end{proof}

\section{Application to $\alpha$-stratified simple modules}\label{s6}

For $c\in\mathbb{C}$ denote by $\Theta_{s}^c=\Theta_{\alpha_s}^c$ 
{\em Mathieu's} twisting functor from \cite[4.3]{Ma}. 

\begin{corollary}\label{ctmain}
Under the assumptions of Theorem~\ref{tmain} we have that the canonical
injection
\begin{displaymath}
U(\mathfrak{g})/\mathrm{Ann}(\Theta_{s}^cL(\mathbf{w}))
\hookrightarrow
\LL(\Theta_{s}^cL(\mathbf{w}),\Theta_{s}^cL(\mathbf{w}))
\end{displaymath}
is surjective. Moreover, for every $w\in W_{\lambda}$ the canonical injection
\begin{displaymath}
U(\mathfrak{g})/\mathrm{Ann}(\Theta_{s}^c\Delta(w))
\hookrightarrow
\LL(\Theta_{s}^c\Delta(w),\Theta_{s}^c\Delta(w))
\end{displaymath}
is surjective.
\end{corollary}

\begin{proof}
From the definition of $\Theta_{s}^c$ it follows that $\Theta_{s}^c$ 
preserves $\mathrm{Ann}(M)$ and induces an isomorphism between
$\LL(M,M)$ and $\LL(\Theta_{s}^c M,\Theta_{s}^c M)$ for any 
$M\in\mathcal{O}$ on which $X_{-\alpha}$ acts injectively.
The first statement now follows from Theorem~\ref{tmain} and the 
second one from Corollary~\ref{cprjan2}.
\end{proof}

When the modules $\Theta_{s}^cL(\mathbf{w})$ are simple, they are simple $\alpha_s$-stratified 
modules considered in \cite{CF,FM}. The modules $\Theta_{s}^c\Delta(w)$ are proper standard
objects in the parabolic generalization of $\mathcal{O}$ studied in
$\cite{FKM}$ (see also \cite{Mz}).

\vspace{1cm}

\begin{center}
\bf Acknowledgments
\end{center}

The research was partially supported by The Royal Swedish Academy of Sciences
and The Swedish Research Council. The research was motivated by some questions, posed to
me by Alexander Stolin, whom I would like to thank very much. I would also like 
to thank Catharina Stroppel for several very stimulating discussions.

\vspace{0.5cm}


\noindent 
Department of Mathematics, Uppsala University, SE-751 06, Uppsala, SWEDEN,
e-mail: {\small \tt mazor@math.uu.se}, 
web: http://www.math.uu.se/$\tilde{\hspace{1mm}}$mazor/
\vspace{0.5cm}



\end{document}